\documentclass[12pt]{iopart}

\usepackage{epsfig,verbatim}

\newcommand{\eps}{\varepsilon}
\newcommand{\ba}{\begin{array}}
\newcommand{\ea}{\end{array}}
\newcommand{\dis}{\displaystyle}
\newcommand{\rens}{\mathbb{R}}
\newcommand{\zens}{\mathbb{Z}}
\newcommand{\torus}{\mathbb{T}}
\newcommand{\Jac}{\mbox{D}}

\usepackage{amssymb}
\usepackage{amsfonts}
\newcommand{\mod}{\mbox{ mod }}

\newtheorem{theo}{\bf Theorem}[section]
\newtheorem{rem}{\bf Remark}[section]
\newenvironment{proof}{\noindent \textit{Proof:}}{\hfill $\square$}

\begin{document}

\title{Non-ergodicity of Nos\'e-Hoover dynamics}
\author{Fr\'ed\'eric Legoll$^{1,2}$, Mitchell Luskin$^3$ and Richard Moeckel$^3$}
\address{
$^1$
Universit\'e Paris-Est,
Institut Navier, LAMI, \'Ecole Nationale des Ponts et Chauss\'ees,
6 et 8 avenue Blaise Pascal, 77455 Marne-La-Vall\'ee Cedex 2, France
}
\address{
$^2$
INRIA Rocquencourt, MICMAC Team-Project, Domaine de Voluceau, B.P. 105,
78153 Le Chesnay Cedex, France
}
\address{
$^3$
School of Mathematics, University of Minnesota, 206 Church Street SE,
Minneapolis, MN 55455, U.S.A.
}
\eads{\mailto{legoll@lami.enpc.fr},
\mailto{luskin@math.umn.edu},
\mailto{rick@math.umn.edu}}

\begin{abstract}
The numerical integration of the Nos\'e-Hoover dynamics
gives a deterministic method that is used to sample the canonical Gibbs measure.
The Nos\'e-Hoover dynamics extends the physical Hamiltonian dynamics by
the addition of a ``thermostat'' variable, that is coupled nonlinearly with
the physical variables. 
The accuracy of the method
depends on the dynamics being ergodic.  Numerical experiments
have been published earlier that are consistent with
non-ergodicity of the dynamics for some model problems.  The authors
recently proved the non-ergodicity of the Nos\'e-Hoover dynamics for the
one-dimensional harmonic oscillator.  

In this paper, this result is extended to
non-harmonic one-dimensional systems.  It is also shown for some
multidimensional systems that the averaged dynamics for the limit
of infinite thermostat "mass" have many invariants, thus giving
theoretical support for either non-ergodicity or slow ergodization.
Numerical experiments for a two-dimensional central force problem
and the one-dimensional pendulum problem give evidence for non-ergodicity.
\end{abstract}





\maketitle

\section{Introduction}
The computation of equilibrium statistical properties of
molecular systems is of great importance in materials science, computational
physics, chemistry, and biology~\cite{frenkelsmit,mcquarrie}.
These equilibrium statistical
properties are given by phase space integrals of the
form
\begin{equation}
\label{ps_aver}
\langle A \rangle = \int A(q,p) \, d\mu(q,p),
\end{equation}
where $q=(q_1,\ldots,q_N) \in {\mathbb R}^{N}$ and $p=(p_1,\ldots,p_N)
\in {\mathbb R}^{N}$ denote a set of positions and
momenta and $A(q,p)$ is an observable, a function defined over
the phase space and related to the macroscopic quantity under
study.
The computation of
integrals such as~(\ref{ps_aver}) is often a challenging problem, especially
when the number of degrees of freedom is large.

For molecular systems at fixed temperature $\theta,$
the measure $d\mu$ is the Gibbs measure for the canonical ensemble
\cite{frenkelsmit,mcquarrie}
\begin{equation}
\label{measure_can}
d\mu(q,p) =
\left[ \frac{\exp\left({-\beta H(q,p)}\right)}
{\displaystyle{ \int \exp\left({-\beta H(q,p)}\right) \ dq \, dp }} \right]
\, dq \, dp,
\end{equation}
where $H(q,p)$ is the Hamiltonian of the system and
$\beta$ is related
to the temperature $\theta$ by $\beta = 1/(k_B \theta)$ with $k_B$ denoting
the Boltzmann constant.
We will consider Hamiltonians of the general form
\begin{equation}
\label{H_qp}
H(q,p) = \frac {p^TM^{-1}(q)p}2+ V(q),
\end{equation}
where $M(q)\in {\mathbb R}^{N\times N}$ for $q\in {\mathbb R}^{N}$
is the generalized mass matrix and $V(q)$ is the potential energy.
We assume that the generalized mass matrix $M(q)\in {\mathbb R}^{N\times N}$
is symmetric and positive definite, so its inverse 
$M^{-1}(q)\in {\mathbb R}^{N\times N}$
exists for all $q\in {\mathbb R}^{N}$ and is also symmetric and positive
definite.

Many methods have been proposed and utilized to approximate the
phase space integral~\eref{ps_aver}, including methods based on
 stochastic or
deterministic dynamics for $(q,p).$
If the dynamics is ergodic with respect to the measure $d\mu$ given by
(\ref{measure_can}), then
the phase-space average
(\ref{ps_aver}) is equal to the time average
\begin{equation}
\label{ergo0}
\int A(q,p) \, d\mu(q,p)
=
\lim_{T \to +\infty} \frac{1}{T} \int_0^T A \left( q(t),p(t) \right) dt
\end{equation}
over a trajectory $(q(t),p(t))_{t \geq 0}$. Thus, the time average can be
approximated by
\begin{equation*}
\lim_{T \to +\infty} \frac{1}{T} \int_0^T A \left( q(t),p(t) \right)
dt \approx \lim_{{\mathcal N} \to +\infty} \frac{1}{\mathcal N}
\sum_{\ell=1}^{\mathcal N}
A(q_\ell,p_\ell),
\end{equation*}
where $(q_\ell,p_\ell)_{\ell \geq 1}$ is a numerical solution of the chosen
dynamics.

In this paper, we investigate the deterministic dynamics
known as Nos\'e-Hoover dynamics~\cite{Hoover}, which is
still widely used although variants have been developed with the goal
to improve its efficiency and overcome its
deficiencies~\cite{Martyna92,Tuckerman00,NPoincare99,RMT05}.  
This dynamics has been first proposed in the form of a
Hamiltonian dynamics on an extended phase space~\cite{nose84}, the
Hamiltonian being chosen such that the marginal distribution of its
microcanonical density is the canonical Gibbs density for the physical
variables. The Nos\'e-Hoover dynamics is then
constructed by rescaling time and momentum to obtain a non-Hamiltonian 
dynamics with physical time and momentum~\cite{Hoover}.

Stochastic
dynamics (such as the Langevin equation, or the recently proposed
Hoover-Langevin method~\cite{LeNoTh09}) can also be considered. See 
\cite{comparisonNVT} for a review
of sampling methods of the canonical ensemble, along with a
theoretical and numerical comparison of their performances for molecular
dynamics.

The equality \eref{ergo0} relies on an ergodicity condition.
This condition has been rigorously proven neither for the
Nos\'e-Hoover dynamics, nor for any other deterministic method
commonly used in practice. In fact, there is numerical evidence that
shows that the Nos\'e-Hoover method is not ergodic for some systems
\cite{Hoover,Martyna92,Tuckerman00}, including the one-dimensional
harmonic oscillator. In \cite{LLM}, we have rigorously analyzed the
dynamics in this special case, and indeed proven the non-ergodicity, for
some regime of parameters.

In this article, we study more general systems. After
briefly recalling the Nos\'e-Hoover equations (see Section
\ref{background}), we first
consider a class of multidimensional systems (see Section
\ref{multiD}). Taking the limit of an infinite thermostat ``mass'' in the
Nos\'e-Hoover equations, we formally obtain an averaged dynamics,
for which we prove the existence of many invariants.
These theoretical results are illustrated by numerical simulations of a
specific system (see Section \ref{multiD_num}). We numerically observe
that, for finite thermostat mass, these invariants are of course
not exactly preserved, but still remain close to
their initial value. This prevents the Nos\'e-Hoover system from
thermalizing. In Section \ref{1Dgeneral},
we next turn to the one-dimensional case, for which we obtain stronger
results. We first prove non-ergodicity of the Nos\'e-Hoover
dynamics, when the mass of the thermostat is large enough (see
Section \ref{sec:1Dtheory}). Our
method extends the one we used to study the harmonic oscillator
case \cite{LLM}. Section \ref{1pendulum} describes an example of such a
one degree of freedom
problem. Again, numerical simulations illustrate the obtained theoretical
results.

\section{Nos\'e-Hoover dynamics}
\label{background}

The Nos\'e-Hoover dynamics involves the physical variables
$q$ and $p$ and
one additional scalar variable, $\xi$, which
represents the momentum of a thermal bath exchanging
energy with the system. The differential equations are:
\begin{equation}
\label{nh-dyn}
\begin{array}{rcl}
\dot{q} &=& \dis{
\frac{\partial H}{\partial p}=M^{-1}(q)p
},
\\
\dot{p} &=& \dis{
-\frac{\partial H}{\partial q}- \frac{\xi}{Q}\, p=-\nabla V(q) -\frac
{p^T\nabla M^{-1}(q)p}2- \frac{\xi}{Q}\, p
},
\\
\dot{\xi} &=& \dis{
p^TM^{-1}(q)p - \frac{N}{\beta}
},
\end{array}
\end{equation}
where $\dot\_$ denotes the time-derivative. The parameter $Q$
represents the mass of the thermostat; it is a free parameter that the
user has to choose.

We recall that invariant measures $\rho(z) \, dz$ for a general
dynamical system
\begin{equation*}
\dot{z}=f(z)
\end{equation*}
are determined by the equilibrium equation
\begin{equation*}
\mbox{div}(\rho(z) f(z))=0.
\end{equation*}
It can be verified by direct computation that the dynamics
\eref{nh-dyn} preserves the measure
\begin{equation}
\label{invmeasure}
d\mu_{\rm NH} = \exp \left[ -\beta \left( H(q,p) + \frac{\xi^2}{2Q}
\right) \right] \ dq \ dp \ d\xi
\end{equation}
by using the fact that the kinetic energy
$\dis \frac{p^TM^{-1}(q)p}2$ is quadratic in $p.$

If the dynamics (\ref{nh-dyn}) is ergodic with respect to $d\mu_{\rm
NH}$, then, by integrating out $\xi$, we have that
the dynamics $(q(t),p(t))$ is ergodic with respect to the
Gibbs measure.  In this case, the time-average of a function $A(q,p)$
along a typical
Nos\'e-Hoover trajectory provides an estimate for
the space-average of $A$ with respect to Gibbs measure.
Unfortunately, the system is generally not ergodic.  In \cite{LLM} we
proved non-ergodicity in the case of the one-dimensional harmonic
oscillator. Our aim here is to study more general
systems.

\section{Systems with first integrals}
\label{multiD}

In this section, we show how the presence of additional integrals for a
Hamiltonian system can impede ergodization of the Nos\'e-Hoover
dynamics.

\subsection{Homogeneous integrals}
\label{sec:homogeneous_integrals}
Consider a Hamiltonian system
\begin{equation}
\label{complete}
\dot{q} =\frac{\partial H}{\partial p},
\quad
\dot{p} = -\frac{\partial H}{\partial q},
\end{equation}
for energy (\ref{H_qp}) which admits a first integral other than $H$
itself.  This means that there is a smooth function $F(q,p)$ whose
Poisson bracket with $H$ vanishes, i.e.,
\begin{equation}
\label{Poissonbracket}
\{H,F\} = H_q^T F_p - H_p^T F_q=0.
\end{equation}
If $F$ is a homogeneous function of the momentum variables, then it gives
rise to a first integral of the Nos\'e-Hoover system.

\begin{theo}\label{thhomogeneous}
If $F(q,p)$ is a first integral of  (\ref{H_qp}) which is homogeneous of
degree $k$ with respect to the momentum variables, $p$, then
\begin{equation}
\label{homogeneous}
G(q,p,\xi) = \frac{\xi^2}{2Q} + H(q,p) - \frac{N}{\beta k}\ln |F(q,p)|
\end{equation}
is a first integral of the corresponding Nos\'e-Hoover system (\ref{nh-dyn}).
\end{theo}

The proof is a simple computation using (\ref{Poissonbracket}),
(\ref{nh-dyn}) and the fact that $F_p^T\,p = k F$.

\smallskip

Of course, the existence of such an integral immediately gives
non-ergodicity of the Nos\'e-Hoover system with respect to
\eref{invmeasure}. For a simple example of a system admitting such a
homogeneous integral, see Section~\ref{multiD_num}.

\subsection{Completely integrable systems and action-angle variables}
\label{sec:action_angle_multiD}

We now assume that the Hamiltonian dynamical system (\ref{complete})
is completely integrable, i.e., the system admits $N$ independent first
integrals which commute in the sense that the Poisson brackets of any
two of them vanish \cite{Arn}. The rest of this section is devoted to
showing that these integrals, even if they are not homogeneous, have a
deleterious effect on the ergodization of the corresponding
Nos\'e-Hoover system.

The non-degenerate level sets of the integrals are $N$-dimensional
manifolds and if they are compact then their connected components  are
diffeomorphic to the $N$-dimensional torus $\torus^{N}$. Moreover, such
a torus has a neighborhood $U\subset \rens^{N}\times \rens^{N}$ in which
one can introduce symplectic action-angle variables. More precisely,
there exist angle variables $\theta\in\torus^{N},$
action variables $a\in D \subset \rens^{N}$, and a symplectic diffeomorphism
$\psi:U\to \torus^{N} \times D$
which transforms (\ref{complete}) to the form
\begin{equation}
\label{eq:newton_action}
\dot{\theta} = \omega(a),
\quad
\dot{a} = 0.
\end{equation}
Here $D$ is an open subset of $\rens^{N}$.
Equivalently, the action-angle Hamiltonian
$\tilde H(\theta,\,a)=H(\psi^{-1}(\theta,\,a))$ is independent of
$\theta$ and
$\dis \frac{\partial \tilde H(\theta,\,a)}{\partial a}=\omega(a).$

In what follows, it will be convenient to define the angle mapping
$\psi_1(q,\,p)\in \rens^{N}$ and the action mapping
$\psi_2(q,\,p)\in\rens^{N}$ by
$\psi(q,\,p)=(\psi_1(q,\,p),\,\psi_2(q,\,p))\in\rens^{2N}.$
We will also use the abbreviated notation
$\dis \partial_1 \psi_i=\frac{\partial \psi_i}{\partial q}
\in\rens^{N\times N}$ and
$\dis \partial_2 \psi_i=\frac{\partial \psi_i}{\partial p}
\in\rens^{N\times N}.$ We then denote the Jacobian of $\psi$ by
\begin{equation*}
\Jac \psi= \left(
\begin{array}{rcl}
\partial_1\psi_1 & \partial_2\psi_1
\\
\partial_1\psi_2 & \partial_2\psi_2
\end{array} \right).
\end{equation*}

We denote the inverse mapping of $\psi(q,p)$
by $\phi(\theta,a).$ The matrix $(\Jac \psi)=(\Jac \phi)^{-1}$
has a simple form since $\phi$ is symplectic:
\begin{equation}\label{sym5}
(\Jac \psi)
=J^{-1} (\Jac \phi)^T J = \left( \ba{cc} (\partial_2
\phi_2)^T & -(\partial_2 \phi_1)^T \\ - (\partial_1 \phi_2)^T &
(\partial_1 \phi_1)^T \ea \right),
\end{equation}
where
\begin{equation*}
J= \left(
\begin{array}{rcl}
0 & I_N
\\
-I_N & 0
\end{array} \right).
\end{equation*}

The diffeomorphism $\psi$ transforms \eref{complete} to
\eref{eq:newton_action} by the chain rule
\begin{equation}
\label{first2}
\left( \ba{c}
\dot{\theta} \\ \dot{a} \ea \right)=
\Jac \psi \  \left( \ba{c} \dot{q} \\ \dot{p} \ea \right).
\end{equation}
Since $\psi$ is symplectic, the dynamics \eref{eq:newton_action} is
obtained from the Hamiltonian $\tilde H(\theta,\,a)$. Hence,
\begin{equation}
\label{first3}
\left( \ba{c}
\omega(a) \\ 0 \ea \right)=
\left( \ba{c} \frac{\partial \tilde H}{\partial a}
\\ -\frac{\partial \tilde H}{\partial \theta} \ea \right)
= \Jac \psi \  \left( \ba{c} \frac{\partial H}{\partial p}
\\ -\frac{\partial H}{\partial q} \ea \right).
\end{equation}

\subsection{Recasting the Nos\'e-Hoover dynamics}

We now multiply the Nos\'e-Hoover equations~\eref{nh-dyn} for $\dot{q}$
and $\dot{p}$ by $\Jac \psi$
to obtain from \eref{first2} and \eref{first3} that
\begin{equation*}
\left( \ba{c} \dot{\theta} \\ \dot{a} \ea \right) =
\left( \ba{c} \omega(a) \\ 0 \ea \right)
- \frac{\xi}{Q} \ \Jac \psi \ \left( \ba{c} 0 \\ p \ea \right).
\end{equation*}
As we are interested in the regime $Q \gg 1$, we rescale by
\begin{equation*}
\eps = 1/\sqrt{Q}, \quad \alpha = \xi/\sqrt{Q}.
\end{equation*}
Using the symplectic property \eref{sym5} of $\psi,$ we
see that the Nos\'e-Hoover equation (\ref{nh-dyn}) can be then
given in the scaled angle-action variables by
\begin{equation}
\label{eq:NH_action_angle}
\left( \ba{c} \dot{\theta} \\ \dot{a} \\ \dot{\alpha} \ea
\right) =
\left( \ba{c} \omega(a) \\ 0 \\ 0 \ea \right)
+ \eps \left( \ba{c} \alpha (\partial_2 \phi_1)^T \phi_2 \\
-\alpha (\partial_1 \phi_1)^T \phi_2 \\
\phi_2^TM^{-1}(\phi_1)\phi_2 - N \beta^{-1} \ea \right).
\end{equation}

\subsection{Averaging the fast variables}

We next apply the averaging method to obtain an approximate system which
does not involve the fast variables. Rigorous results about averaging
for Hamiltonian systems with several degrees of freedom are fraught with
technical difficulties (see for example~\cite{AKN}). These arise from
the fact that $\omega(a)$, the frequency vector of the fast angles,
experiences resonances of the form $k\cdot \omega(a)=0$,
$k\in\zens^N$, for certain values of the action vector $a$. In fact,
these resonant actions are generally dense in the action domain $D$. In
spite of this, the averaged differential equations often provide a
useful first approximation to the behavior of the slow variables when
$\varepsilon$ is small.

In our problem, the averaged system for the slow variables is given by
\begin{equation}
\label{averaged}
\begin{array}{rcl}
\dot{a} &=& - \alpha \ S(a),
\\
\dot{\alpha} &=& k(a),
\end{array}
\end{equation}
where
\begin{equation}
\label{eq:def_S_k}
\begin{array}{rcl}
S(a) &=&\dis{
\langle (\partial_1 \phi_1)^T \ \phi_2 \rangle(a)
=
\int_{\torus^{N}}(\partial_1 \phi_1)^T (\theta,\,a)
\phi_2(\theta,\,a)\,d\theta
},
\\
k(a) &=& \dis{
\langle \phi_2^TM^{-1}(\phi_1)\phi_2 \rangle(a) -
\frac{N}{\beta}
}
\\
&=& \dis{
\int_{\torus^{N}} \phi_2^T(\theta,\,a) M^{-1}
(\phi_1(\theta,\,a))\phi_2(\theta,\,a)
\,d\theta- \frac{N}{\beta} 
}.
\end{array}
\end{equation}

We next show that, with an additional assumption on the action-angle
mapping $\phi$, we have
\begin{equation}\label{aa}
S(a) = a.
\end{equation}
Recall that a map $(q,p) = \phi(\theta,a)$ is symplectic if it preserves
the canonical differential two-forms, i.e.,
$\phi^*(\sum dp_i\wedge dq_i) = \sum da_i\wedge d\theta_i$
where $\phi^*$ indicates the pull-back, meaning that we write
$p_i, q_i, dp_i, dq_i$ in terms of the $\theta, a$ variables. It follows
from this that the difference of the corresponding canonical one-forms
$\phi^*(p^T\,dq) - a^T\,d\theta$ is a closed one-form on $ \torus^{N}
\times D$. We will call $\phi$ {\em exact symplectic} if this closed
one-form is exact, i.e., if
\begin{equation}\label{exsym}
\phi^*(p^T\,dq) = a^T\,d\theta +dF(\theta, a)
\end{equation}
where $F(\theta,a)$ is a real-valued function on $\torus^{N}\times D$. This
stronger condition holds for the action-angle coordinates associated to
many well-known integrable systems.

As an example, consider the one-dimensional harmonic oscillator $H(q,p)=
\frac12(p^2+q^2)$ which is a completely integrable system with $N=1$
degrees of freedom. In any annulus of the form $0<h_1\le (p^2+q^2)/2 \le
h_2$, we can introduce action-angle variables $(\theta,a)$ such that
\begin{equation*}
q = \sqrt{2a}\cos\theta, \qquad p = -\sqrt{2a}\sin\theta.
\end{equation*}
For $\phi(\theta,a) = (q,p)$, we have
\begin{equation*}
\phi^*(p\,dq)= -\sin\theta \cos\theta da +2a \sin^2\theta d\theta
\end{equation*}
and so
\begin{equation*}
\phi^*(pdq) - ad\theta =  -\sin\theta \cos\theta da +a(2\sin^2\theta-1)
d\theta = dF
\end{equation*}
where
\begin{equation*}
F= -a\sin\theta\cos\theta.
\end{equation*}

It turns out that the action-angle variables constructed according to
the usual method of Arnold \cite{Arn} always have this exactness property.
To see this, recall that in Arnold's method, the tori given by fixing
the $N$-independent integrals of motion are parametrized by
angle variables $\theta= (\theta_1,\ldots,\theta_N)$ derived from the
$N$ commuting Hamiltonian flows defined by the integrals. Then the
action variables are given by
\begin{equation}\label{actionvars}
a_i = \int_{\gamma_i} p^T\,dq
\end{equation}
where $\gamma_i$ is the curve in the torus defined by holding $\theta_j
= \mbox{const}$ for $j\ne
i$ and letting $\theta_i$ run over $[0,1]$.  The integral depends
on which torus is
considered, i.e., it is a function of the $N$ first integrals.
The usual proof shows that the map
$(q,p)=\phi(\theta,a)$ is defined and symplectic on some domain of the form
$\torus^{N}\times D$. It follows that
\begin{equation*}
\nu  = \phi^*(p^T\,dq) - a\,d\theta
\end{equation*}
is a closed differential one-form on $\torus^{N}\times D$. Showing that
$\phi$ is exact
symplectic amounts to showing that $\nu$ is exact. For this, it suffices
to check that
its integral around any closed curve vanishes. In fact, since $\nu$ is
closed, it suffices to check the curves of the form
$C_i = \Gamma_i \times \{a_0\}, a_0\in D$
where $\Gamma_i = \{\theta:\theta_j= \mbox{const}, j\ne i\}$.
For such a curve, we have
\begin{equation*}
\int_{C_i} \nu = \int_{C_i} \phi^*(p\,dq) - a\,d\theta =  \int_{\gamma_i}  p\,dq -
\int_0^1 a_{0i}\,d\theta_i =  a_{0i}-a_{0i} = 0.
\end{equation*}
Here we used the fact that under $\phi$ the curve $C_i$ maps to the  curve $\gamma_i$
used in (\ref{actionvars}).

To prove \eref{aa} under the exact symplectic assumption  \eref{exsym}, we first note that
\begin{equation*}
\phi^*(p^T\,dq) = \phi_2^T d\phi_1 =
\phi_2^T(\partial_1\phi_1\,d\theta + \partial_2\phi_1\,da).
\end{equation*}
Hence, \eref{exsym} reads
\begin{equation*}
\phi_2^T(\partial_1\phi_1\,d\theta + \partial_2\phi_1\,da) =
a^T\,d\theta +dF(\theta, a).
\end{equation*}
For any $j=1,\ldots,N$, we can integrate both sides with respect to
$\theta_j$, along the circular
loop $C_j$ as in the last paragraph. Using the
periodicity of $F$ in $\theta_j$, we obtain
\begin{equation}\label{explicit}
\sum_{i=1}^N\int_0^{1}\left(\frac {\partial \phi_{1i}}{\partial \theta_j}\right)
(\theta,\,a) \phi_{2i}(\theta,\,a)\,d\theta_j=a_j.
\end{equation}
We can then further integrate \eref{explicit} over the angles
$\theta_k$ for $k\ne j$ to obtain that
\begin{equation*}
S(a)_j=
\sum_{i=1}^N\int_{\torus^{N}}\left(\frac {\partial \phi_{1i}}{\partial \theta_j}\right)
(\theta,\,a) \phi_{2i}(\theta,\,a)\,d\theta=a_j
\end{equation*}
for $j=1,\dots,N.$

\subsection{First integrals of the averaged Nos\'e-Hoover equations}
\label{sec:1st_int}

A direct calculation shows that a set of $N$ independent first integrals
for the averaged Nos\'e-Hoover equations
\begin{equation}\label{goodform}
\begin{array}{rcl}
\dot a &=& -\alpha a,\\
\dot \alpha &=& k(a),
\end{array}
\end{equation}
are given by
\begin{equation}\label{integ}
\begin{array}{rcl}
G_i(a,\,\alpha) &=& \dis{ \frac{a_i}{a_N},\qquad i=1,\dots,N-1, }
\\
G_N(a,\,\alpha)&=& \dis{
\frac{\alpha^2}2+\int^{a_N}\frac{k\left(s\frac{a_1}{a_N},
\dots,s\frac{a_{N-1}}{a_N},s\right)}s\,ds}.
\end{array}
\end{equation}
To prove that $G_N(a,\,\alpha)$ is a first integral, it is helpful to use
the fact that $G_i(a,\,\alpha)=a_i/a_N$ for $i=1,\dots,N-1$ is a first integral.

We summarize the result of this section in the following theorem.
\begin{theo} \label{thaveraged}
The averaged equations for the Nos\'e-Hoover dynamics for a
completely integrable Hamiltonian system
has $N$ independent first integrals.
\end{theo}

To the extent that the averaging method applies, we expect that $G_i(a(t),\,\alpha(t))$
evolves slowly for small $\eps$ and so the sampling of the Gibbs measure is slow
even if the dynamics is ergodic. We will verify this numerically in an example in the
next section.
It turns out that $G_i(a(t),\,\alpha(t))$ remains quite close to its
initial value for fairly large values of $\eps$ as well.

\section{A central force problem}
\label{multiD_num}

We consider here a two degrees of freedom system to illustrate the
theoretical results obtained in the previous section. We work with the
Hamiltonian (\ref{H_qp}) with $N=2$, the identity mass matrix $M(q) =
I_2$, and a potential $V(|q|)$ which depends only on the distance to the origin.
The Hamiltonian system (\ref{complete}) admits two first integrals, the
energy $H$ and the angular momentum
\begin{equation*}
L = q_1 p_2 - q_2 p_1,
\end{equation*}
which satisfy $\{H,L\} = 0$, and whose gradients are linearly
independent, except for values of $H$ and $L$ satisfying a condition of
the form $f(H,L) = 0$ for some function $f$. Hence, this system is
completely integrable.
Assume
that $V(r) \rightarrow +\infty$ as $r \rightarrow +\infty$. Then level
sets of $H$ are compact, hence the level sets
$\{ (q,p) \in \rens^4; \ H(q,p) = h, L(q,p) = \ell \}$ are also compact,
hence there
exists action-angle variables for this system.

To describe the action variables, first introduce polar coordinates
$(r,\phi)$ in $\rens^2$. The angular momentum is
\begin{equation*}
L= r^2 \dot\phi.
\end{equation*}
Fixing a value for $L$, we have a reduced Hamiltonian system for the
radial variables $(r,p_r)$, where $p_r=\dot r$, with Hamiltonian
\begin{equation*}
H_L(r,p_r) = \frac12 p_r^2 + \frac12 \frac{L^2}{r^2} + V(r).
\end{equation*}
This reduced system has one degree of freedom and can be understood by
the usual phase-plane method. Since $V(r)
\rightarrow +\infty$ as $r \rightarrow +\infty$, the level curves
\begin{equation*}
C_{(h,L)} = \{(r,p_r):H_L(r,p_r) = h\}, \qquad L\ne 0,
\end{equation*}
generically consist of one or more simple closed curves. For the
unreduced system, where we remember the angle $\phi$, each such curve
becomes an invariant torus $T_{(h,L)}$. It can be shown that the action
variables assigned to such a torus by Arnold's procedure are as follows:
$a_1(h,L)$ is the area in the $(r,p_r)$ plane enclosed by the simple
closed curve of the reduced system,
\begin{equation*}
a_1(h,L) = \int_{C_{(h,L)}} p_r \, dr,
\end{equation*}
and $a_2(h,L) = L$, the angular momentum. Note that $a_1(h,L)$ is easily
computable by standard numerical integration schemes.

Since $L$ is homogeneous of degree $k=1$ in the momentum variables,
Theorem~\ref{thhomogeneous} gives a first integral for the Nos\'e-Hoover
system:
\begin{equation}
\label{eq:inv1}
G(q,p,\xi) = \frac{\xi^2}{2Q} + H(q,p) -
\frac{2}{\beta} \ln \left| L(q,p) \right|.
\end{equation}
In addition, Theorem~\ref{thaveraged} provides additional integrals for
the averaged Nos\'e-Hoover equations, in particular the ratio of the
action variables
\begin{equation}
\label{eq:inv2}
G_1(q,p) = \frac{a_1(q,p)}{a_2(q,p)} = \frac{a_1(H(q,p),L(q,p))}{L(q,p)}.
\end{equation}
To the extent that the averaging method applies for this two-degrees of
freedom problem, this ratio should evolve only very slowly when $Q \gg
1$. In the sequel, we present some
numerical simulations showing first that, when $Q \gg 1$, this is indeed
the case, and second that, for $Q =1$, such a behaviour persists
to some extent.

\medskip

\begin{figure}[htbp]
\centerline{\input{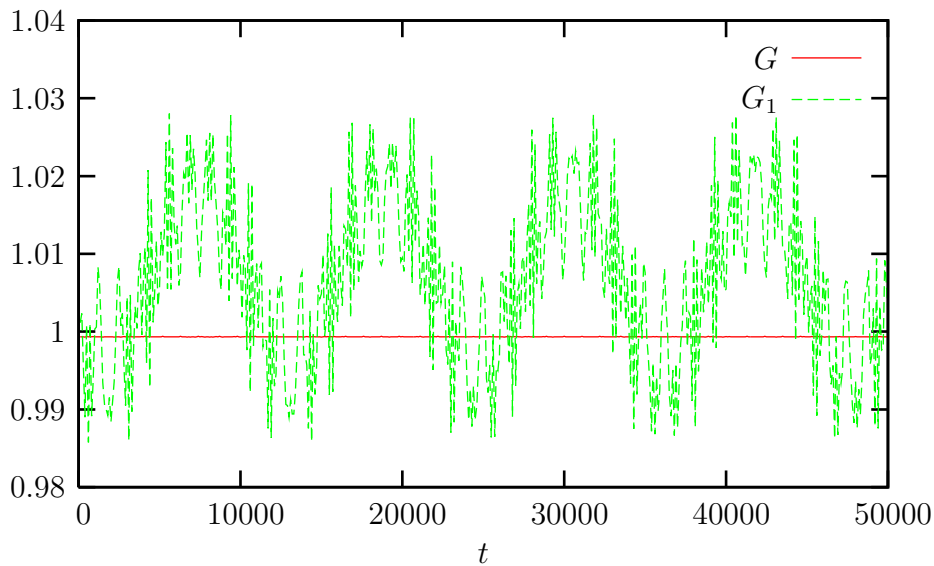}}
\caption{Plot of $G(q(t),p(t),\xi(t))$ and $G_1(q(t),p(t))$
 (renormalized by their initial value) along the trajectory of
(\ref{nh-dyn}), for $Q=100$ ($\beta = 1$, initial
condition $q=(0;0.5)$, $p=(-1.5; 1.5)$, $\xi=0$).}
\label{fig:Q100}
\end{figure}

Consider the example with potential
\begin{equation*}
V(r) = r^2+r^4,
\end{equation*}
with an initial condition $(q_0,p_0,\xi_0)$ such that
$L(q_0,p_0) \neq 0$. We compute the trajectory of the Nos\'e-Hoover
dynamics \eref{nh-dyn} with the algorithm proposed in \cite{molphys96}.
On Figure~\ref{fig:Q100}, we plot
$G(q(t),p(t),\xi(t))$ and $G_1(q(t),p(t))$, where $G$ and $G_1$ are
defined by \eref{eq:inv1} and \eref{eq:inv2}, for $Q=100$. We indeed
observe that $G$ is preserved, whereas $G_1$ evolves slowly.

\begin{figure}[htbp]
\centerline{\input{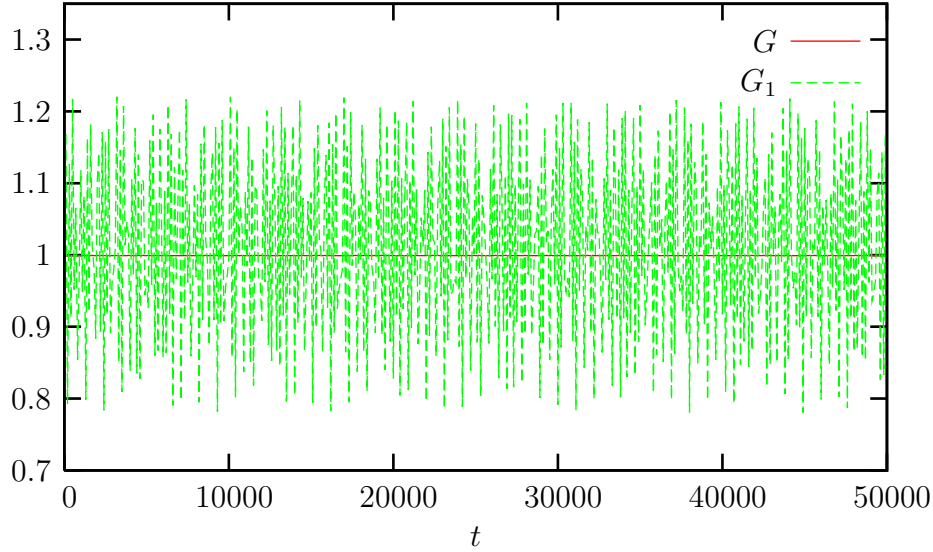}}
\caption{Plot of $G(q(t),p(t),\xi(t))$ and $G_1(q(t),p(t))$
 (renormalized by their initial value) along the trajectory of
(\ref{nh-dyn}), for $Q=1$ ($\beta = 1$, initial
condition $q=(0;0.5)$, $p=(-1.5; 1.5)$, $\xi=0$).}
\label{fig:Q1}
\end{figure}

We now consider the value $Q=1$ and plot the same quantities as above on
Figure~\ref{fig:Q1}. Again, $G$ is preserved, whereas $G_1$ evolves in a
band which is still quite narrow, even for this small value of $Q$.

\medskip

Let us now derive another
quantity, which does not behave as well as $G_1$
for large $Q$, but happens to behave in a better way for small
$Q$\footnote{We have a clear
  understanding of why  $G_1$ behaves better than this quantity when $Q
  \gg 1$. However, the situation for $Q=1$ is less clear.}.
From the Nos\'e-Hoover dynamics \eref{nh-dyn}, we compute that
\begin{eqnarray}
\label{Ldyn}
\dot{L} &= - \frac{\xi}{Q} L = - \eps \alpha L,
\\
\label{Hdyn}
\dot{H} &= \dis{ - \frac{\xi}{Q} p^T p = - \eps \alpha \, p^T p
= - \eps \alpha \, \phi_2^T(\theta,a) \, \phi_2(\theta,a) }.
\end{eqnarray}
Since $\theta$ are fast variables whereas $a$, $\alpha$ and $H$ are slow
ones, we again formally use the averaging method on \eref{Hdyn} and consider the dynamics
\begin{equation}
\label{eq:H_approx}
\dot{H} = - \eps \alpha \, k_0(a)
\end{equation}
with
\begin{equation*}
k_0(a)= \int_{\torus^2} \phi_2^T(\theta,a) \, \phi_2(\theta,a) \,
d\theta
\end{equation*}
(note that \eref{Ldyn} does not depend on the fast variables $\theta$).
Now recall that the action variables $a$ are functions of $H$ and
$L$. The equation \eref{eq:H_approx} hence reads
\begin{equation*}
\dot{H} = - \eps \alpha \, k_0(H,L).
\end{equation*}
The averaged system is thus
\begin{equation}
\label{averaged_2d}
\begin{array}{rcl}
\dot{L} &=& - \alpha L,
\\
\dot{H} &=& - \alpha \, k_0(H,L),
\\
\dot{\alpha} &=& k_0(H,L) - 2 \beta^{-1}.
\end{array}
\end{equation}
Note that
\begin{equation}
\label{eq:def_E}
E(H,L,\alpha) = \frac{\alpha^2}2 + H -
\frac{2}{\beta} \ln \left| L \right|
\end{equation}
is a first integral of the above system. It is just the analogue of
\eref{eq:inv1} being a first integral for the Nos\'e-Hoover system (see
Theorem~\ref{thhomogeneous}).

\begin{figure}[htbp]
\centerline{\input{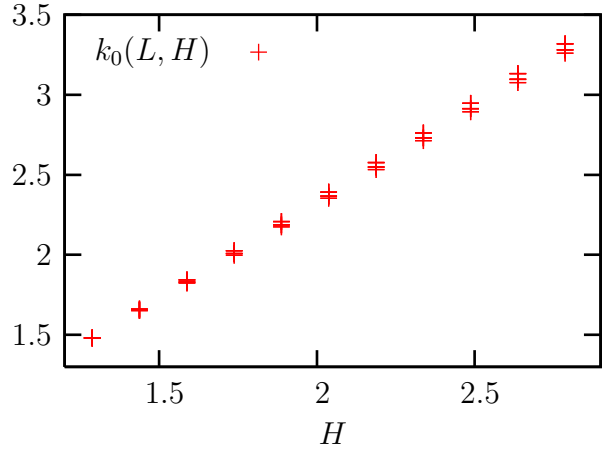}}
\caption{Plot of $H \mapsto k_0(L,H)$. For each value of $H$, we have
 considered several values of $L$.}
\label{fig:k0_2d}
\end{figure}

On Figure~\ref{fig:k0_2d}, we plot the function $H \mapsto k_0(H,L)$,
for several values of $L$. We observe that $k_0(H,L)$
is almost a constant with respect to $L$, and can hence be
approximated\footnote{In practice, we have considered several energy values
  $H_i$, and for each $H_i$, we have considered several configurations
  $(q_{i,j},p_{i,j})$ with energy $H_i$ and angular momentum
  $L_{i,j}$. We next have computed $k_0(H_i,L_{i,j})$ by averaging
  $p(t)^T p(t)$ along a constant energy trajectory. Averaging these
  $k_0(H_i,L_{i,j})$, we obtain $k_0^{\rm app}(H_i)$, which next leads
  to $k_0^{\rm app}(H)$ for any $H$ by piecewise linear interpolation.}
by a function $k_0^{\rm app}(H)$.  

We hence approximate \eref{averaged_2d} by
\begin{equation}
\label{averaged_2d_bis}
\begin{array}{rcl}
\dot{L} &=& - \alpha L,
\\
\dot{H} &=& - \alpha \, k_0^{\rm app}(H),
\\
\dot{\alpha} &=& k_0^{\rm app}(H) - 2 \beta^{-1}.
\end{array}
\end{equation}
Now, it is natural to introduce the variable $\tau$ defined by
\begin{equation*}
\tau(H) = \exp \left( \int^{H} \frac{ds}{k_0^{\rm app}(s)} \right)
\end{equation*}
and its reciprocal $H(\tau)$, such that \eref{averaged_2d_bis} reads
\begin{equation}
\label{averaged_2d_ter}
\begin{array}{rcl}
\dot{L} &=& - \alpha L,
\\
\dot{\tau} &=& - \alpha \, \tau,
\\
\dot{\alpha} &=& k_0^{\rm app}(H(\tau)) - 2 \beta^{-1}.
\end{array}
\end{equation}
This system is in the form \eref{goodform}. Its two first integrals
are
\begin{equation*}
E_1(L,\tau) = \frac{\tau}{L}
\end{equation*}
and
\begin{eqnarray*}
E_2(\tau,\alpha) &=&
\frac{\alpha^2}2 + \int^\tau \frac{k_0^{\rm app}(H(s)) - 2
 \beta^{-1}}{s} \, ds
\\
&=&
\frac{\alpha^2}2 + \int^\tau \frac{k_0^{\rm app}(H(s))}{s} \, ds -
\frac{2}{\beta} \ln \tau
\\
&=&
\frac{\alpha^2}2 + H(\tau) - \frac{2}{\beta} \ln \tau.
\end{eqnarray*}
The first invariant $E$ given by \eref{eq:def_E}
is not independent from $E_1$ and $E_2$: $E_2 = E - 2 \beta^{-1} \ln
|E_1|$.

\begin{figure}[htbp]
\centerline{\input{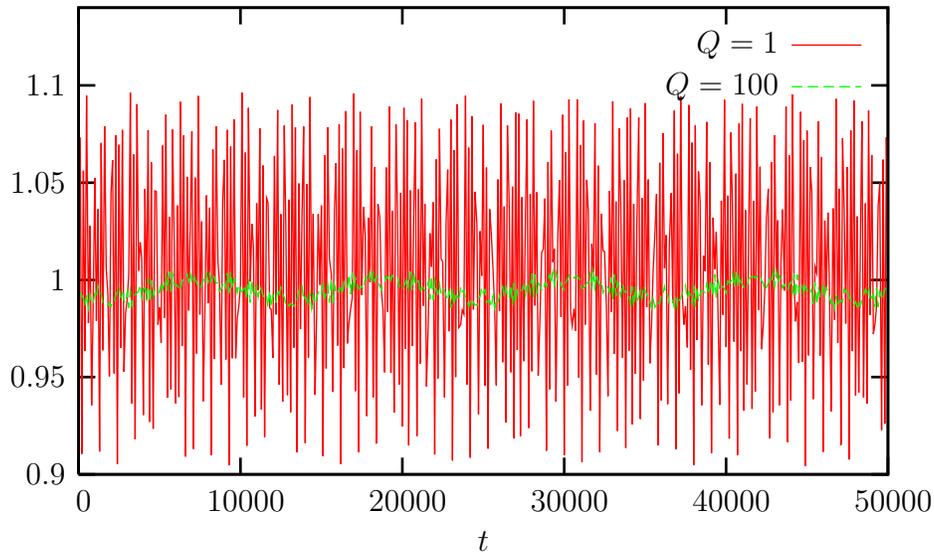}}
\caption{$E_1(q(t),p(t))$ (renormalized by its initial
  value) along the trajectory of
 (\ref{nh-dyn}), for $Q=1$ and $Q=100$ ($\beta = 1$, initial
 condition $q=(0;0.5)$, $p=(-1.5; 1.5)$, $\xi=0$).}
\label{fig:G1_particular}
\end{figure}

We now consider the same trajectories of the Nos\'e-Hoover dynamics
that we considered on Figures~\ref{fig:Q100} and \ref{fig:Q1}, and we
plot
\begin{equation*}
E_1(q,p) = E_1 \left( L(q,p),\tau(H(q,p))\right).
\end{equation*}
We see on Figure~\ref{fig:G1_particular} that, for $Q=100$, this
quantity is almost preserved, and that, even for $Q=1$, it remains close
to its initial value.

\begin{figure}[htbp]
\centerline{\input{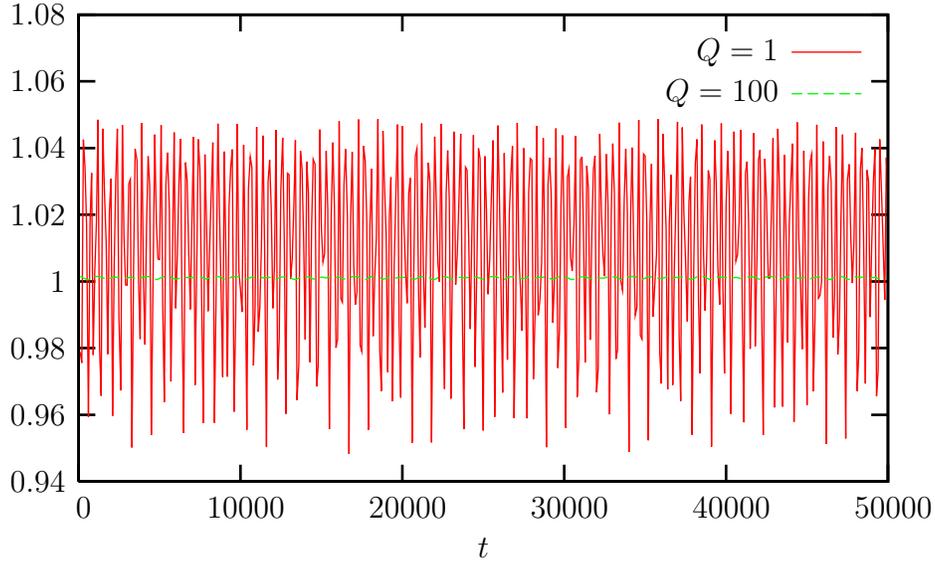}}
\caption{$E_2(q(t),p(t),\xi(t))$ (renormalized by its initial
  value) along the trajectory of
 (\ref{nh-dyn}), for $Q=1$ and $Q=100$ ($\beta = 1$, initial
 condition $q=(-0.5;0.5)$, $p=(-1; 1)$, $\xi=0$).}
\label{fig:G2_particular}
\end{figure}

We finally consider an initial condition such that $L(q_0,p_0) =
0$. Along the trajectory of \eref{nh-dyn}, we have
$L(q(t),p(t))=0$ by \eref{Ldyn}, hence $E_1$ is not defined.
On Figure~\ref{fig:G2_particular}, we plot
\begin{equation*}
E_2(q,p,\xi) = E_2 \left( \tau(H(q,p)), \frac{\xi}{\sqrt{Q}} \right)
\end{equation*}
along two trajectories, obtained with the same initial condition and the
choices $Q=100$ and $Q=1$. We again observe that $E_2$ is almost
constant for $Q=100$, and that it remains close to its
initial value for $Q=1$.

\begin{figure}[htbp]
\centerline{\input{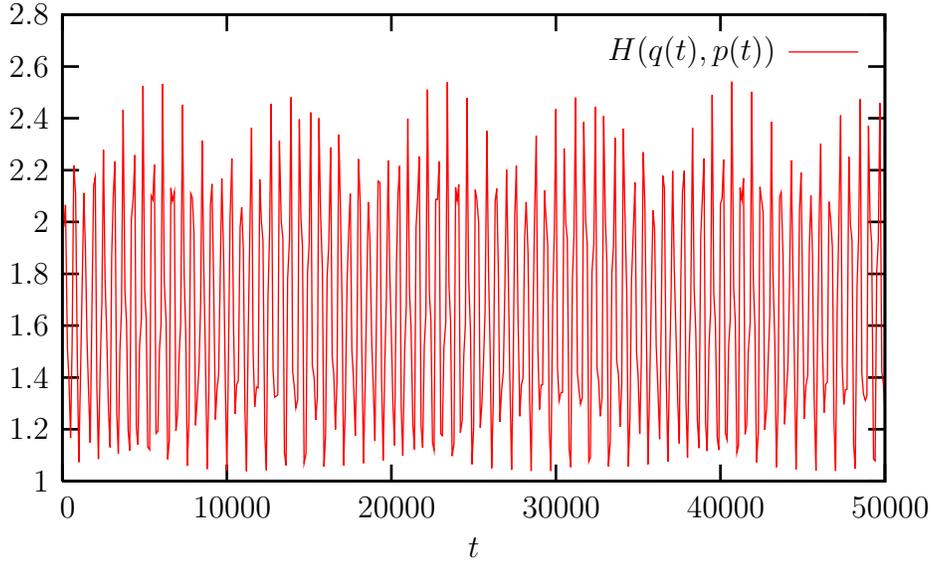}}
\caption{$H(q(t),p(t))$ along the trajectory of
(\ref{nh-dyn}), for $Q=1$ ($\beta = 1$, initial
condition $q=(-0.5;0.5)$, $p=(-1; 1)$, $\xi=0$).}
\label{fig:H_L_null}
\end{figure}

On Figure~\ref{fig:H_L_null}, we plot the energy $H(q(t),p(t))$
along the same trajectory (for $Q=1$). We see that values $h \leq 1$ are
not sampled. However, there exist $(q,p) \in \rens^4$
such that $L(q,p) = 0$ and $H(q,p)$ is as close to 0 as wanted. Hence,
the trajectory only samples a strict subset of the level set
$\{ (q,p); \ L(q,p) = 0 \}$.

\section{Systems with one degree of freedom}
\label{1Dgeneral}
Consider a Hamiltonian system of the form (\ref{H_qp}) with $N=1$ and
$M(q) = 1$.  All such systems are
completely integrable since $H$ itself provides the required integral of motion.  Suppose
there is
an interval of energies $I = [h_1, h_2]$ such that the level curves
$M(h) = \{(q,p):H(q,p)= h\}, h\in I,$ are all simple closed curves (one-dimensional tori)
which are non-degenerate in the sense that
the gradient of $H$ does not vanish. Then the
plane region $U = \{(q,p): h_1\le H(p,q) \le h_2 \}$ is diffeomorphic to an annulus,
and we can introduce action-angle variables $(a,\theta)$ in $U$, and an
exact symplectic map $\phi(\theta,a) = (q,p)$, as in Section
\ref{sec:action_angle_multiD}.

From Theorem~\ref{thaveraged}, we have $N=1$ first integrals for the averaged
Nos\'e-Hoover equations. Let us now rewrite
this first integral more explicitly. In view of \eref{integ}, we have
\begin{equation}
\label{eq:G}
G(a,\alpha) = \frac{\alpha^2}{2} + W(a)
\end{equation}
with
\begin{equation}
\label{eq:W1D}
W(a) = \int^a \frac{k(s)}{s} \, ds,
\end{equation}
where, in view of \eref{eq:def_S_k}, $k$ is given by
\begin{equation}
\label{eq:k1D}
k(a) = \int_{\torus}
\phi^2_2(\theta,\,a) \, d\theta - \frac{1}{\beta}.
\end{equation}

As in the multidimensional case, this integral prevents rapid
ergodization, at least for small $\eps$.  But in the one-degree of
freedom case we go further and identify conditions on $H$ which
rigorously imply non-ergodicity. The method is essentially the one
used in \cite{LLM} where we treated the harmonic oscillator.  Namely,
the integral $G$ leads to invariant tori of the averaged system which,
under certain assumptions, persist for small values of $\eps$.

\subsection{Proof of non-ergodicity}
\label{sec:1Dtheory}

We will apply a KAM theorem to the Nos\'e-Hoover equations, in the formulation
\eref{eq:NH_action_angle}. Let us introduce the Poincar\'e return map,
$P_\eps$, of the system \eref{eq:NH_action_angle} to the Poincar\'e
section defined by $\theta = 0 \mod 1$. It is convenient to rescale time
by $\omega(a)$, so that the return time when $\eps = 0$ is 1. This just
alters the parametrization
of the solutions so that the return time to the Poincar\'e section is~1.

Since there is only one degree of freedom, the averaging method can be rigorously
justified.
Indeed we can eliminate the fast angle
$\theta$ of (\ref{eq:NH_action_angle}) by a change of variables.
We construct functions $g(\hat{a},\theta,\hat{\alpha})$
and $h(\hat{a},\theta,\hat{\alpha})$ and corresponding new variables
$(\hat{a},\hat{\alpha})$ defined by
\begin{equation}
\label{eq:change0}
\begin{array}{rcl}
a &=& \hat{a} + \eps g(\hat{a},\theta,\hat{\alpha}),
\\
\alpha &=& \hat{\alpha} + \eps h(\hat{a},\theta,\hat{\alpha}),
\end{array}
\end{equation}
so that in
the new variables $(\hat{a},\hat{\alpha})$, the dynamics
(\ref{eq:NH_action_angle}) is given (after replacing
$(\hat{a},\hat{\alpha})$ by $({a},{\alpha})$) by
\begin{equation}
\label{eq:NH_action_angle20}
\begin{array}{rcl}
\dot{\theta}&=&\omega(a)+O(\eps),\\
\dot{a}&=&-\eps \alpha S(a)+O(\eps^2),\\
\dot{\alpha}&=&\eps k(a)+O(\eps^2),
\end{array}
\end{equation}
where $S(a)$ and $k(a)$ are the averages (\ref{eq:def_S_k}).

In view of \eref{eq:change0} and \eref{eq:NH_action_angle20},
the Poincar\'e map $P_\eps(\hat{a},\hat{\alpha})$ is an $O(\eps^2)$
perturbation of the time $\eps$ advance map of the
averaged system \eref{averaged}, for which $G$ defined
by \eref{eq:G} is a first integral.
So we now make some assumptions about
the level curves of $G$.

Recall that we are working in a region $U$ of the
$(q,p)$-plane defined by an interval of energies $I$ which corresponds to an interval of
actions
$J = [a_1,a_2] $. We assume that $W(a)$ has at least one local minimizer
$a_0$ in $J$. We have
\begin{equation*}
0 = W'(a_0) = \frac{k(a_0)}{a_0},
\end{equation*}
hence $k(a_0) = 0$ and the point $P=(a_0,0)$ is an equilibrium point for
(\ref{averaged}). The parts of the level curves of $G$ which are near
$P$ are simple closed curves around $P$ in the $(a,\alpha)$-plane.

\begin{rem}
\label{rem:eq2solve}
If $a_0$ is a local minimizer of $W$, then $k(a_0) = 0$, hence
$\dis \int_{\torus} \phi^2_2(\theta,\,a_0) \, d\theta  = \beta^{-1}$.
\end{rem}

\begin{figure}[htbp]
\centerline{\input{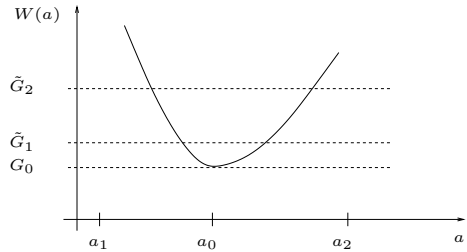}}
\caption{Schematic representation of the interval $\tilde{K} = [\tilde
G_1,\tilde G_2]$.}
\label{fig:tilde_K}
\end{figure}

Let $G_0 = G(a_0,0) = W(a_0)$ be the value of the integral $G$ at the
equilibrium point $P=(a_0,0)$. Choose constants $\tilde{G_1}$ and
$\tilde{G_2}$ such that $G_0 < \tilde{G_1} < \tilde{G_2} <
\min(G(a_1,0),G(a_2,0)) = \min(W(a_1),W(a_2))$
and let $\tilde K = [\tilde G_1,\tilde G_2]$ (see
Figure~\ref{fig:tilde_K}). Then the level curves $\{
(a,\alpha); \ G(a,\alpha) = c \}$, where $c \in \tilde{K}$, have connected
components which are simple closed curves near $P$. The union of these
components for $c \in
\tilde{K}$ forms a region $\tilde D$ near $P$ which is diffeomorphic to
an annulus.

The variable $G$ defines a natural action variable in $\tilde D$.
We construct the corresponding angle variable $\phi$ by following the
same method as in \cite{LLM}. Let $T_1(g)$ denote
the period of the periodic solutions of (\ref{averaged}) which
corresponds to the level curve $G = g \in \tilde K$. The averaged
differential equation (\ref{averaged}) becomes
\begin{equation}
\label{averaged3}
\begin{array}{rcl}
\dot{\phi} &=& 1/T_1(G),
\\
\dot{G} &=& 0.
\end{array}
\end{equation}
In these coordinates, the time $\eps$ advance map takes the form
$(\phi,G) \mapsto (\phi_1,G_1)$ where
\begin{equation}
\label{timeeps}
\begin{array}{rcl}
\phi_1 &=& \phi + \eps/T_1(G),
\\
G_1 &=& G.
\end{array}
\end{equation}
Call this map $Q_\eps(\phi,G)$. Then $P_\eps(\phi,G) = Q_\eps(\phi,G) +
O(\eps^2)$.

\begin{theo}
\label{theo:non_ergo}
Suppose the period function $T_1(G)$ is not identically
constant on the interval $\tilde K$. Then, for $\eps$ sufficiently
small, the Poincar\'e map $P_\eps$ has invariant circles in the region
$\tilde D$ and so the Nos\'e-Hoover system is not ergodic: trajectories
$(q(t),p(t))$ that solve (\ref{nh-dyn}) are not ergodic with respect to
the Gibbs measure \eref{measure_can}.
\end{theo}

\begin{proof}
As in \cite{LLM}, we apply Moser's twist theorem to the Poincar\'e
map $P_\eps$. The details are similar to those in \cite{LLM}, so we
only sketch the argument here.

The fact that the Nos\'e-Hoover differential equation preserves the
invariant measure \eref{invmeasure} implies (as in \cite{LLM}) that
$P_\eps$ preserves an
invariant measure in the $(a,\alpha)$ plane. It follows that the maps
$P_\eps$ have the curve intersection property. The hypothesis on
$T_1(G)$ guarantees that, making $\tilde K$ and hence $\tilde D$ smaller
if necessary, we may assume that either $T_1'(G)>0$ or $T_1'(G)<0$
throughout $\tilde D$. This means that there are many invariant circles
in $\tilde D$ for which the rotation number under $Q_\eps$ is Diophantine.
Moreover, the required twist condition holds. Moser's theorem guarantees
that such invariant circles perturb to nearby invariant circles for
$\eps$ sufficiently small.

We now show that the existence of these invariant circles implies
non-ergodicity with respect to the Gibbs measure. First, note that
on a level curve ${\mathcal M} = \left\{ (a,\alpha); \ G(a,\alpha) = c
\right\}$, where $c \in \tilde{K}$, we have that $W(a)$ is bounded from
above, since $W(a) \leq G(a,\alpha) = c$. In view of the choice of
$\tilde{K}$ (see Figure~\ref{fig:tilde_K}), this implies that $a$ is
lower and upper bounded. Since $a'(h)$ is positive and bounded away from
0 for $h \in I = [h_1,h_2]$, this hence shows that $H(q,p)$ is lower and upper bounded
(that is,
$| H(q,p)|$ is bounded) on the invariant circle ${\mathcal M}$ of
$Q_\eps$.
As a consequence, $|
H(q,p)|$ is bounded on the nearby invariant circles of $P_\eps$. Hence,
the trajectory of (\ref{nh-dyn}) does not sample values of $H(q,p)$
larger than some threshold. This is a contradiction with $(q(t),p(t))$
that solves (\ref{nh-dyn}) being ergodic with respect to the Gibbs
measure \eref{measure_can}.
\end{proof}

\medskip

Because of the complicated series of coordinate changes leading from the
original Hamiltonian system to the averaged system, it is not easy to
state simple conditions on the original potential function $V(q)$ which
guarantee that the period function $T_1(G)$ is not constant.    An
equilibrium point surrounded by periodic orbits of constant period is
called isochronous and various criteria for isochronicity have been
given.  Our problem can be reduced to a Hamiltonian case for which a
simple criterion can be stated.

To carry out the reduction, replace $a$ in (\ref{averaged}) (that is,
\eref{goodform}) by $\sigma = \ln (a/a_0)$, where $a_0$ is a local
minimizer of $W$ (see Figure~\ref{fig:tilde_K}). The
differential equation (\ref{goodform}) becomes
\begin{equation}
\label{averaged4}
\dot \sigma = -\alpha,
\quad
\dot \alpha = U'(\sigma),
\end{equation}
where $U(\sigma) = W(a_0 \exp \sigma)$. Except for a reversal of time, this
is a classical Hamiltonian system with Hamiltonian
$G(\sigma,\alpha) = \alpha^2/2 + U(\sigma)$. It has an
equilibrium point at the origin $(\sigma,\alpha) = (0,0)$.

Now \cite{LL} discusses the problem of recovering the potential of such
a system from its period function (see also \cite{CMV}).  Let $G_0 =
U(0)$ be the energy level of the equilibrium point at the origin.  For
$G>G_0$ let $L(G)$ be the width of the potential well at energy $G$,
i.e., $L(G)= \sigma_2(G)-\sigma_1(G)$ where $\sigma_i(G)$ are the two
roots of $U(\sigma) = G$ near $\sigma=0$. Then $T_1(G)$ is constant
if and only if $\displaystyle L(G) = \frac{T_1}{\pi}\sqrt{2(G-G_0)}$.
This is just the formula for the width of the quadratic potential well
associated to a harmonic oscillator of period $T_1$. Clearly this is highly
exceptional and is easy to rule out, at least numerically.

Another way to show that $T_1(G)$ is non-constant is to observe that the
constancy of the period implies that the family of periodic orbits
surrounding the equilibrium point must fill the entire plane \cite{CMV}. If
this were not so, then there would be another equilibrium point on the
boundary of the maximal family which would force
$T_1(G)\rightarrow\infty$. For example, for certain values of
$\beta$, the pendulum equations (see Section \ref{1pendulum}) lead
to an averaged system with more than one equilibrium, and this
immediately implies that $T_1(G)$ is non-constant.

\subsection{The simple pendulum problem}
\label{1pendulum}

\begin{figure}[htbp]
\centerline{\input{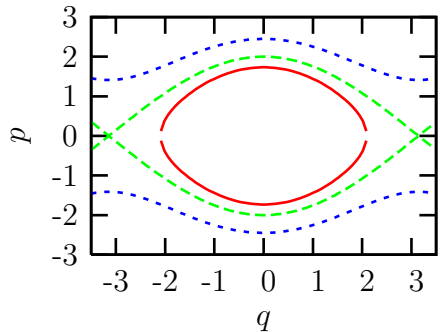}}
\caption{Phase portrait of the simple pendulum.}
\label{fig:pendulum}
\end{figure}

We consider here the numerical example of a simple pendulum whose
potential energy is given by
\begin{equation*}
V(q) = - \cos q.
\end{equation*}
We reduce $q$ modulo $2 \pi$.
By construction, the energy satisfies $h \geq -1$. The phase portrait is
shown on Figure~\ref{fig:pendulum}.
The above assumptions are satisfied for energies
in the interval $I=[h_1,h_2]$, with $-1<h_1<h_2<1$, or $1<h_1<h_2$.

\begin{figure}[htbp]
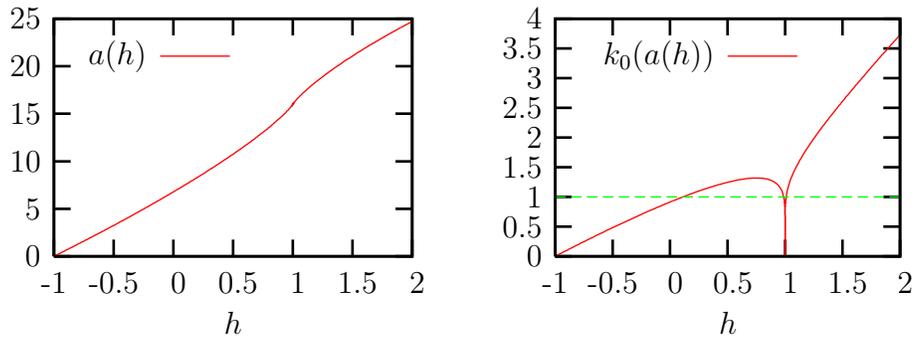

\centerline{
\input{figure9a.tex}
\input{figure9b.tex}
}
\caption{Numerically computed values of $a(h)$ and $k_0(a(h))$ (see text).}
\label{fig:a_h}
\end{figure}

First, we numerically compute $a(h)$ defined by
\eref{actionvars}. Note that, in this one-dimensional setting,
\begin{equation*}
a(h) = \int_{M(h)} p \ dq ,
\end{equation*}
where the line integral is taken in the direction of the Hamiltonian
flow, and $M(h) = \{(q,p) \in \rens^2: \ H(q,p)= h\}$. We also compute
\begin{equation*}
k_0(a) = \int_{\torus} \phi^2_2(\theta,\,a) \, d\theta,
\end{equation*}
which is independent of $\beta$ and satisfies $k(a) = k_0(a) - \beta^{-1}$,
with $k$ defined by \eref{eq:k1D}. In practice, $k_0$ is computed using
the fact that, for any energy level $h$,
\begin{equation*}
k_0(a(h)) = \lim_{T \to +\infty} \frac{1}{T} \int_0^T p^2(t) \, dt,
\end{equation*}
where $(q(t),p(t))$ solve the Newton equations of motion for the
pendulum at the constant energy $h$. 
Results are shown on Figure~\ref{fig:a_h}.

We have seen that the action values $a$ such that $k_0(a) = \beta^{-1}$,
that is $k(a) = 0$,
play an important role (see Remark \ref{rem:eq2solve}).
In view of Figure~\ref{fig:a_h}, we see that, when $\beta^{-1} < \beta_c^{-1}$ for
some threshold $\beta_c$, then the equation $k_0(a) = \beta^{-1}$ has
three solutions. When $\beta^{-1} > \beta_c^{-1}$, then the equation
$k_0(a) = \beta^{-1}$ has a unique solution. In what follows, we detail
the numerical results obtained with the choice $\beta = 1$, which
corresponds to the first case. Similar results have been
obtained for
choices of $\beta$ corresponding to the other case. Hence, the
conclusions that we draw here are by no means restricted to the case
$\beta = 1$. 

\begin{figure}[htbp]
\centerline{\input{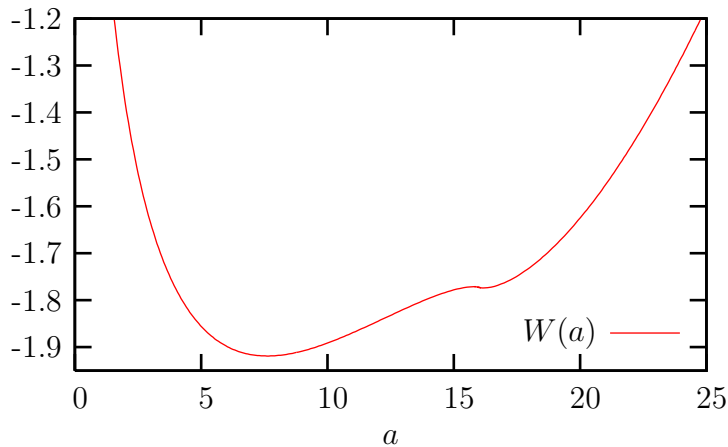}}
\caption{Numerically computed values of $W(a)$ for $\beta = 1$.}
\label{fig:w_h}
\end{figure}

The function $W(a)$ defined by \eref{eq:W1D} is shown on
Figure~\ref{fig:w_h} for the choice $\beta = 1$. This function has two
local minimizers, $a_1 \approx 7.6$ and $a_2 \approx 16.17$. Following
the above theoretical analysis, we work close to one of them. We have
chosen to work close to $a_1$.

\begin{figure}[htbp]
\centerline{\input{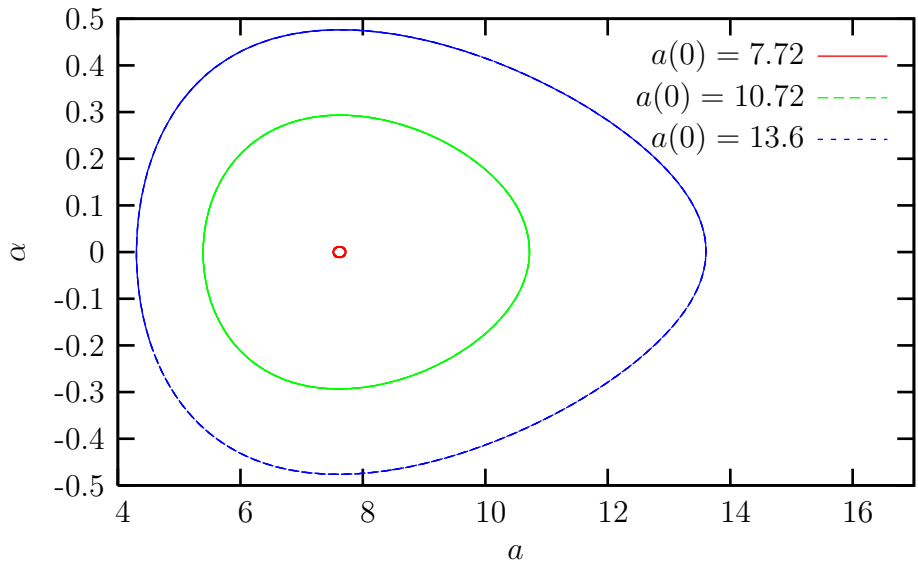}}
\caption{Trajectories of (\ref{averaged}) for several
different values of $a(0)$ ($\beta = 1$).}
\label{fig:level}
\end{figure}

On Figure~\ref{fig:level}, we plot the trajectory of the averaged dynamics
(\ref{averaged}) for different initial conditions. These trajectories have been
computed with the Symplectic Euler algorithm used on the Hamiltonian
formulation (\ref{averaged4}). As expected, the trajectory
is a simple closed curve around the equilibrium point $(a_1,0)$,
that corresponds to a level curve of $G$. These curves are also
invariant curves of the map $Q_\eps$ defined in the previous section
(see map \eref{timeeps}).

\begin{figure}[htbp]
\centerline{\input{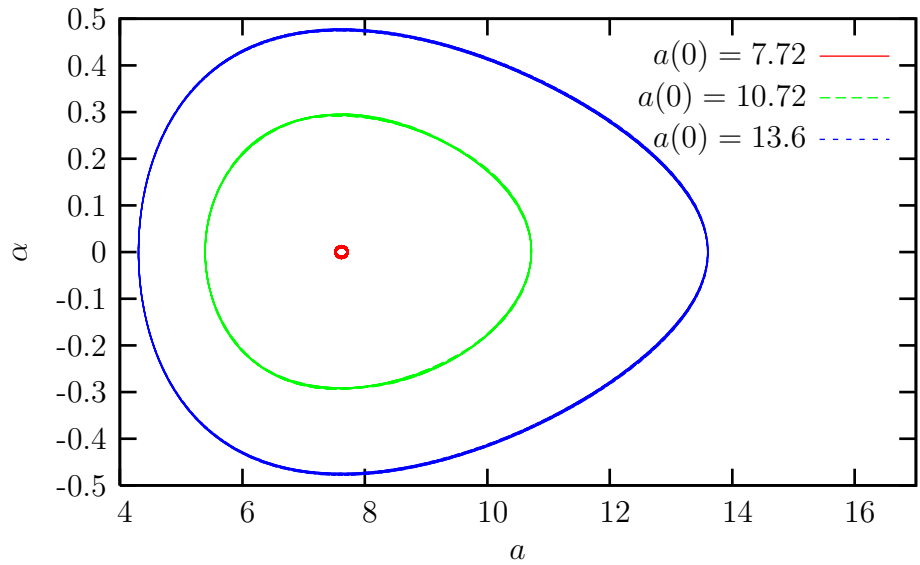}}
\caption{Poincar\'e return map of (\ref{eq:NH_action_angle}) on the plane
$\theta = 0 \mod 1$ for several initial conditions ($Q = 10^5$, $\beta
= 1$).}
\label{fig:level_full}
\end{figure}

\begin{figure}[htbp]
\centerline{\input{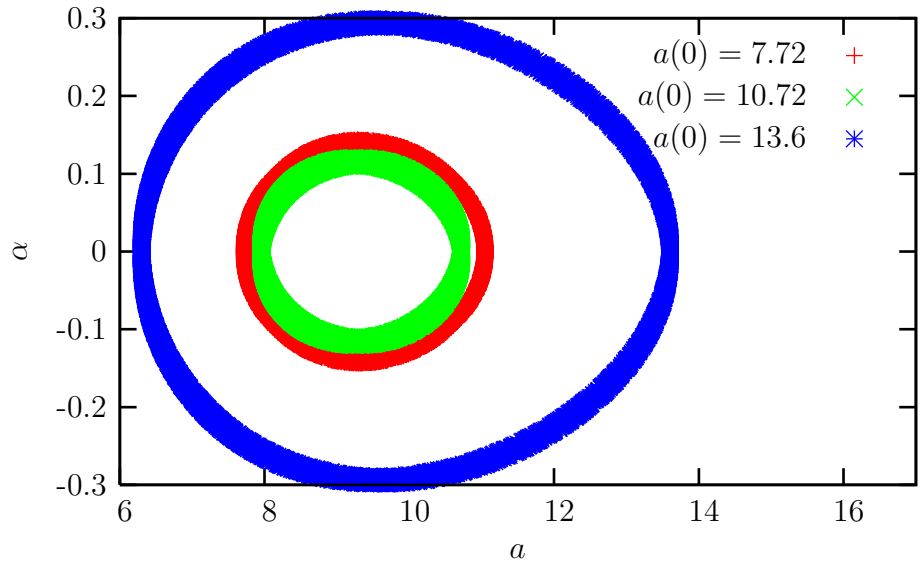}}
\caption{Poincar\'e return map of (\ref{eq:NH_action_angle}) on the plane
$\theta = 0 \mod 1$ for several initial conditions ($Q = 1$, $\beta
= 1$).}
\label{fig:level_full_very_low}
\end{figure}

We now study how these curves persist upon perturbation.
We recall that the Poincar\'e return map $P_\eps$ of the Nos\'e-Hoover
dynamics (\ref{eq:NH_action_angle}) on the section $\theta = 0 \mod
1$ is a perturbation of $Q_\eps$.
Results for the Poincar\'e return map of the dynamics
(\ref{eq:NH_action_angle}) are shown on Figure~\ref{fig:level_full} for $Q
= 10^5$ (that is, $\eps = \sqrt{10} \times 10^{-3}$), and on
Figure~\ref{fig:level_full_very_low} for $Q=1$ (that is, $\eps = 1$),
for the same initial energies as for Figure~\ref{fig:level}.
These Poincar\'e return maps have been computed using the fact that the
section $\theta = 0 \mod 1$ corresponds to the section $q = 0 \mod 2
\pi$.
We see a good agreement between Figures~\ref{fig:level}
and~\ref{fig:level_full}.
The presence of invariant circles on Figures~\ref{fig:level_full} and
\ref{fig:level_full_very_low} shows that the system
(\ref{eq:NH_action_angle}) seems to have invariant curves, for $Q=10^5$
and $Q=1$.

\begin{figure}[htbp]
\centerline{\input{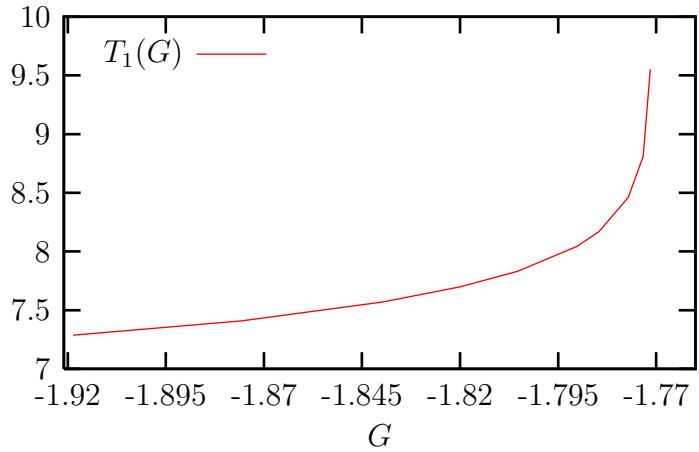}}
\caption{Period $T_1(G)$ of the averaged equation \eref{averaged}
($\beta = 1$).}
\label{fig:t1_g}
\end{figure}

Note that Theorem \ref{theo:non_ergo}, which states the non-ergodicity
of the Nos\'e-Hoover equations, relies on the important assumption
that the period $T_1(G)$ of the averaged equations is not constant. This
holds true for the pendulum case, in view of the discussion at the end of Section
\ref{sec:1Dtheory}. This is
also confirmed by numerical computations of $T_1(G)$ (see
Figure~\ref{fig:t1_g}).

\begin{figure}[htbp]
\centerline{\input{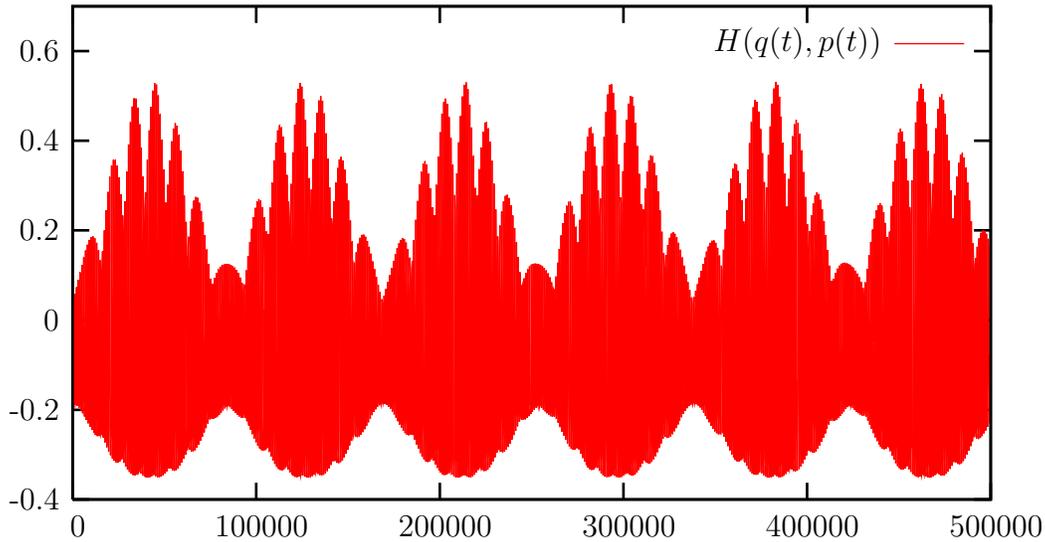}}
\caption{Energy $H(q(t),p(t))$ along the trajectory of
(\ref{nh-dyn}), for $Q=1$ and $\beta = 1$, and the initial
condition $q=0$, $p=1.5$, $\xi=0$ (that is, $a(0) = 7.72$).}
\label{fig:energy_particular}
\end{figure}

Let us now look at another criterion for ergodicity, namely what
energy values are sampled. We see on
Figure~\ref{fig:level_full_very_low} that small values of $a$
are not sampled: we have $a \geq 6$ for the three initial conditions
that we considered. In view of Figure~\ref{fig:a_h}, this
corresponds to small values of $H$ not being sampled. On
Figure~\ref{fig:energy_particular}, we plot the physical energy
$H(q(t),p(t))$ along the trajectory of
(\ref{eq:NH_action_angle}), for the value $Q=1$, and the initial
condition $q=0$, $p=1.5$, $\xi=0$, that corresponds to the initial value
$a(0) = 7.72$ that we studied on Figures~\ref{fig:level},
\ref{fig:level_full} and \ref{fig:level_full_very_low} (results are the same
for other initial conditions). We see that $H \geq
-0.4$. If the dynamics (\ref{eq:NH_action_angle}) was sampling the
canonical measure, then all values of $H$
would be attained. In particular, the smallest values $H \approx -1$ would be
the most frequent ones. Indeed, from the Gibbs measure
\eref{measure_can}, we compute the probability distribution function of
the energy, which reads $\rho(h) = z^{-1} \exp(-\beta h) \, a'(h)$, where
$z$ is a normalization constant. For the pendulum case, $a'(h)$ is
close to a constant (see Figure~\ref{fig:a_h}), hence the smallest
values of $h$ are the most frequent ones.
Hence, it seems that (\ref{eq:NH_action_angle}) is
not ergodic with respect to the canonical measure, even for the
value $Q=1$.

\ack

Part of this work was completed while the first author
was visiting the Institute for Mathematics and its Applications
(Minneapolis), whose hospitality is gratefully acknowledged.

The work of Fr\'ed\'eric Legoll was supported in part by the Agence Nationale de la
Recherche (INGEMOL non-thematic program) and by the Action Concert\'ee Incitative
``Nouvelles Interfaces des Math\'ematiques'' SIMUMOL (Minist\`ere de
  la Recherche et des Nouvelles Technologies, France).

The work of Mitchell Luskin was supported in part by NSF Grants
DMS-0757355 and DMS-0811039, the Institute for Mathematics
and its Applications,
and by the University of Minnesota Supercomputing Institute.
This work is also based on
work supported by the Department of Energy under Award Number
DE-FG02-05ER25706.

Richard Moeckel was partially supported by NSF Grant
DMS-0500443.

\section*{References}

\bibliographystyle{plain}
\bibliography{legoll_luskin_moeckel}

\begin{thebibliography}{10}

\bibitem{Arn}
V.I. Arnold.
\newblock {\em Mathematical Methods of Classical Mechanics}.
\newblock Springer Verlag, New York, 1989.

\bibitem{AKN}
V.I. Arnold, V.V. Kozlov, and A.I. Neishtadt.
\newblock {\em Mathematical Aspects of Classical and Celestial Mechanics}.
\newblock Springer, 2nd ed., 1997.

\bibitem{NPoincare99}
S.D. Bond, B.J. Leimkuhler, and B.B. Laird.
\newblock The {N}os\'e-{P}oincar\'e method for constant temperature molecular
  dynamics.
\newblock {\em J. Comput. Phys.}, 151:114--134, 1999.

\bibitem{comparisonNVT}
E.~Canc\`es, F.~Legoll, and G.~Stoltz.
\newblock Theoretical and numerical comparison of some sampling methods for
  molecular dynamics.
\newblock {\em Math. Mod. Num. Anal. (M2AN)}, 41(2):351--389, 2007.

\bibitem{CMV}
A.~Cima, F.~Manosas, and J.~Villadelprat.
\newblock Isochronicity for several classes of {H}amiltonian systems.
\newblock {\em Jour. Diff. Eq.}, 157:373--413, 1999.

\bibitem{frenkelsmit}
D.~Frenkel and B.~Smit.
\newblock {\em Understanding Molecular Simulation, from algorithms to
  applications, 2nd ed.}
\newblock Academic Press, 2002.

\bibitem{Hoover}
W.~Hoover.
\newblock Canonical dynamics: Equilibrium phase-space distributions.
\newblock {\em Phys. Rev. A}, 31(3):1695--1697, 1985.

\bibitem{LL}
L.D. Landau and E.M. Lifshitz.
\newblock {\em Mechanics}.
\newblock Pergamon Press, Oxford, 2nd ed., 1969.

\bibitem{LLM}
F.~Legoll, M.~Luskin, and R.~Moeckel.
\newblock Non-ergodicity of the {N}os\'e-{H}oover thermostatted harmonic
  oscillator.
\newblock {\em Arch. Rat. Mech. Anal.}, 184(3):449--463, 2007.

\bibitem{LeNoTh09}
B.J. Leimkuhler, E.~Noorizadeh, and F.~Theil.
\newblock private communication.

\bibitem{RMT05}
B.J. Leimkuhler and C.R. Sweet.
\newblock A {H}amiltonian formulation for recursive multiple thermostats in a
  common timescale.
\newblock {\em SIAM J. App. Dyn. Sys.}, 4(1):187--216, 2005.

\bibitem{Martyna92}
G.~Martyna, M.~Klein, and M.~Tuckerman.
\newblock {N}os\'e-{H}oover chains: The canonical ensemble via continuous
  dynamics.
\newblock {\em J. Chem. Phys.}, 97(4):2635--2643, 1992.

\bibitem{molphys96}
G.~Martyna, M.~Tuckerman, D.~Tobias, and M.~Klein.
\newblock Explicit reversible integrators for extended systems dynamics.
\newblock {\em Mol. Phys.}, 87:1117--1157, 1996.

\bibitem{mcquarrie}
D.~McQuarrie.
\newblock {\em Statistical Mechanics}.
\newblock University Science Books, 2000.

\bibitem{nose84}
S.~Nos\'e.
\newblock A unified formulation of the constant temperature molecular dynamics
  method.
\newblock {\em J. Chem. Phys.}, 81(1):511--519, 1984.

\bibitem{Tuckerman00}
M.~Tuckerman and G.~Martyna.
\newblock Understanding modern molecular dynamics: Techniques and applications.
\newblock {\em J. Phys. Chem. B}, 104:159--178, 2000.

\end{thebibliography}

\end{document}